\documentclass{article}

\usepackage{amsthm,amsmath,amssymb}

\begin{document}

\title{A Splitting Criterion for Galois Representations Associated to Exceptional Modular Forms {\mbox{(mod $p$)}}}
\author{Ken McMurdy}
\maketitle

\newcommand{\Qbar}{\overline{\mathbb{Q}}}\newcommand{\Q}{\mathbb{Q}}
\newcommand{\Qpbar}{\overline{\mathbb{Q}}_p}\newcommand{\Qp}{\mathbb{Q}_p}
\newcommand{\Zp}{\mathbb{Z}_p}
\newcommand{\Z}{\mathbb{Z}}
\newcommand{\Fp}{\mathbb{F}_p}\newcommand{\Fpbar}{\overline{\mathbb{F}}_p}
\newcommand{\dia}[1]{{\langle#1\rangle}}

\theoremstyle{plain}
\newtheorem*{conjecture}{Conjecture}
\newtheorem{theorem}{Theorem}[section]
\newtheorem{corollary}{Corollary}[theorem]
\newtheorem{lemma}{Lemma}[theorem]
\theoremstyle{definition}
\newtheorem*{definition}{Definition}
\newtheorem*{note}{Note}
\newtheorem*{remark}{Remark}

\section{Introduction}
\indent
Let $f=\sum a_n q^n$ be a newform of type $(k,\epsilon)$ for $\Gamma_1(N)$ defined over a field $E$ of characteristic $p$ with $(p,N)=1$.
Then there is a construction due to Deligne \cite{De} which attaches to $f$ a continuous semi-simple Galois representation 
$$\rho_f:Gal(\Qbar/\Q)\rightarrow GL_2(E).$$
The representation is unramified outside of $Np$ where for any Frobenius element $\sigma_l$ the characteristic polynomial of $\rho_f(\sigma_l)$ is simply $x^2-a_l x+\epsilon(l)l^{k-1}$.
At the special prime $p$, however, the representation {\em can} be ramified, and the behavior of $\rho_f$ at this prime has been a topic of great interest for many years.

In the ordinary case ($a_p\neq0$) and when $2\leq k\leq p+1$, Deligne has shown that the restriction to a decomposition group at $p$ has the form 
$$\rho_{f,p}=\left(\begin{matrix}
        \chi^{k-1}\cdot\lambda(\epsilon(p)/a_p) & \ast \\
        0                &      \lambda(a_p)
        \end{matrix}\right).$$
Here $\chi$ is the cyclotomic character and $\lambda(\alpha)$ is the unramified character which takes $\sigma_p$ to $\alpha$.
By definition, a Galois representation is unramified if it is trivial on the inertia group and tamely ramified if it is at least trivial on the wild inertia group.
Since $\chi$ is of exponent $p-1$ on the inertia group, we see immediately that $\rho_{f,p}$ can at best be tamely ramified if $k\neq p$.
When $k=p$, it is easily seen that $\rho_{f,p}$ is tamely ramified if and only if it is unramified.
In 1979, Serre conjectured the following purely modular criterion for when these things happen.
\begin{conjecture}[Serre] Suppose $f$ is ordinary ($a_p\neq0$) and $2\leq k\leq p$. The representation $\rho_{f,p}$ will be tamely ramified ($k\neq p$) or unramified ($k=p$) iff there exists an eigenform form $g=\sum b_nq^n$ of type $(k',\epsilon)$ satifying $\theta g=\theta^{k'}f$ for $k'=p+1-k$ (called a ``companion form'').
\end{conjecture}

In 1990, this conjecture was proven by B. Gross \cite{Gr} except in the case where the characters along the diagonal of $\rho_{f,p}$ are not distinct.
In other words, the ``exceptional case'' where $k=p$ and $\epsilon(p)=a_p^2$ was left unproven.
Gross also proved that in these non-exceptional cases the existence of a companion form for $f$ is actually equivalent to the representation being split, i.e. the sum of two characters with respect to some basis.
In 1992 Coleman and Voloch \cite{CV} gave a different proof of the conjecture when $2<k\leq p$.
This proof did not depend on certain unproven compatibilities between cohomology theories which Gross had needed (see the introduction to \cite{Gr}), and it did {\em not} exclude the exceptional case.
However, in the exceptional case it is no longer clear what the relationship is between the ramification of $\rho_{f,p}$ and whether or not it is split.
While the exceptional case with $p=2$ remained an open problem even for the conjecture, in all exceptional cases the question of when $\rho_{f,p}$ is split remained open.

The goal of this paper will be to begin to answer the question of when $\rho_{f,p}$ is split in the exceptional case.
In particular, we will prove a splitting criterion in the case that $H_m$, the completion of the Hecke algebra at the maximal ideal $m=m_f$, is actually equal to $\Zp$.
In terms of modular forms, this condition is equivalent to saying that $f$ has a unique lifting to a newform $F$ of weight $2$ and level $Np$, and that this newform has coefficients in $\Zp$.
One reason for this assumption is that it makes it possible to define a $p$-divisible subgroup $G=G_f$ of $J_1(Np)$ for which $\rho_{f,p}$ is just a twist of the representation on $G[p]$.
As in \cite{Gr} we will attach a splitting invariant $q_p$ to $G[p]$ as well as a characteristic $0$ refinement $q$ to $G$.
The triviality of $q_p$ will determine whether $\rho_{f,p}$ is split, but one very important result of the paper will show how to explicitly calculate $q_p$ from $q$ in the cases under consideration.
The proof of the main theorem then boils down to two formulas for approximating $q$.
The first is the same formula for $d\log q$ which was used in \cite{CV}, involving the cup product on $H^1_{DR}(I_1(N))$.
The second is a formula for $\log q$ involving an inner product on a subspace of $H^0(X_1(Np),\Omega)$ introduced in \cite{C1}.
This inner product can also be interpreted as the cup product on the De Rham cohomology of an algebraic curve which reduces to $I_1(N)$.
With this cohomological interpretation, the main result of this paper is equivalent to the assertion following \cite[Prop. 5.5]{C1}, that when $\rho_{f,p}$ is unramified its splitting is determined by the vanishing of the class ${\overline{[\pi^{-p}F|w]}}$ in $H^1_{DR}(I_1(N))$.

\section{$\rho_{f}$ and the $p$-divisible Group $G$} 
In this section we will recall how the newform $f$ gives rise to a Galois representation $\rho_f$ as well as a $p$-divisible group $G$.
As always we begin with an ordinary newform $f$ of type $(k,\epsilon)$ for $\Gamma_1(N)$ over a finite field $E$ of characteristic $p$ with $(N,p)=1$.
Furthermore, we assume $p>2$.
Then there is a lifting of $f$ to a newform $F$ of type $(2,\epsilon_F)$ for $\Gamma_1(Np)$ over the integral closure $R$ of $\Zp$ in $\Qpbar$ \cite[Prop. 9.3]{Gr}.
\begin{note}If $F(q)=\sum A_n q^n$ and $f(q)=\sum a_n q^n$ this means that $A_n\equiv a_n\pmod{m_R}$.
Also, if $\epsilon_F=\epsilon_N\cdot\epsilon_p$ we must have $\epsilon_N(d)\equiv\epsilon(d)\pmod{m_R}$.
By necessity the guaranteed lifting must also satisfy $\epsilon_p=(\chi_{{}_T})^{k-2}$, where $\chi_{{}_T}:(\Z/p\Z)^*\rightarrow \Zp^*$ is the Teichm\"{u}ller character.\end{note}

Now, let $K$ be the finite extension of $\Qp$ over which $F$ is defined.
Let $H$ be the commutative subring of $End(J_1(Np))$ generated by the Hecke operators $T_l$ for $l$ prime to $Np$ and $\dia{d}_{Np}$.
By the multiplicity one theorem \cite[Thm. 6.2.3]{DI} and the fact that $F$ is a newform, the subspace $W$ of $T_p(J_1(Np))\otimes K$ on which $H\otimes K$ acts via the character associated to $F$ is a $2$-dimensional $K$-vector space.
Therefore from $W$ we obtain a representation
$$r_F:Gal(\Qbar/\Q)\rightarrow GL_2(K)$$
One can then use the Eichler-Shimura congruence, $T_l\equiv l/\sigma_l+\dia{l}_{Np}\sigma_l$, to calculate the characteristic polynomial of $r_F$ at a Frobenius element $\sigma_l$ with $l\nmid Np$.
We define $\rho_F=r_F\otimes\epsilon_F$ so that the characteristic polynomial of $\rho_F(\sigma_l)$ is simply $x^2-A_lx+l\epsilon_F(l)$.
Finally, since $\rho_F$ stabilizes an $\mathcal{O}_K$-lattice we may reduce modulo $m_K$ and define $\rho_f$ to be the semisimplification of this reduction.
It is then easy to show that the characteristic polynomial of $\rho_f(\sigma_l)$ is $x^2-a_lx+\epsilon(l)l^{k-1}$ for $l$ prime to $Np$.

\begin{note}
It is {\em almost} possible to construct the representation $\rho_f$ without first lifting $f$.
Let $V_f$ be the subspace of $J_1(Np)[p](\Qbar)\otimes E$ on which $T_l$ acts like $a_l$ for $l$ prime to $Np$, $\dia{d}_N$ acts like $\epsilon(d)$, and $\dia{d}_p$ acts like $d^{k-2}$.
\cite[Prop. 11.8]{Gr} states that when $\rho_f$ is irreducible, the semi-simplification of $V_f\otimes\epsilon\chi^{k-2}$ is isomorphic to a direct sum of copies of $\rho_f$.
The problem is that $f$ can have multiple liftings.
So even if we consider the full Hecke algebra, $V_f$ can still have dimension greater than $2$.
Again this is one of the reasons why we will make the assumption that $H_m=\Zp$ which implies that the lifting $F$ is unique.
\end{note}

For the definition of $G$ then we begin by enlarging $H$ to include also the operators $U_l$ for $l|Np$.
Following \cite{Gr} there is a maximal ideal $m$ of $H$ with $H/m=E$ such that modulo $m$ we have the following congruences.
\begin{align*}
    T_l&\equiv a_l\\
    U_l&\equiv a_l\\
    <d>_N&\equiv \epsilon(d)\\
    <d>_p&\equiv d^{k-2}
\end{align*}
Let $H_m$ and $H_p$ be the completions of $H$ at $m$ and $p$, and let $\epsilon_m$ be the idempotent of $H_p$ satisfying $H_m=\epsilon_m H_p$.

The Tate module of $J_1(Np)$ is a module for $H_p$ and the submodule $\epsilon_m T_p(J_1(Np))$ is both free over $\Zp$ and stable under $Gal(\Qbar/\Q)$.
Therefore it defines a $p$-divisible group $G$ over $\Q$ with $T_p(G)=\epsilon_m T_p(J_1(Np))$ which is acted on by $H_m$.
In fact, in many cases $T_p(G)$ is free over $H_m$ which gives $G$ the structure of an $m$-divisible group.
At this point, however, we will make the assumption that $H_m=\Zp$ so that the two notions are equivalent.
The following theorem, almost a direct quote of \cite[Prop. 12.9]{Gr}, tells us what the Galois structures of $G$ and $G[m]=G[p]$ are in that case.
\begin{theorem}\label{Gstr}In the case where $H_m=\Zp$, the $p$-divisible group $G$ satisfies the following:\\
1. $G$ has height $2$ and good reduction over $\Zp[\zeta_p]$. In the canonical exact sequence $0\rightarrow G^0\rightarrow G\rightarrow G^e\rightarrow0$ over $\Zp[\zeta_p]$, $G^0$ and $G^e$ both have height $1$.\\
2. The filtration $0\rightarrow T_p(G^0)\rightarrow T_p(G)\rightarrow T_p(G^e)\rightarrow0$ is stable under all of $\mathcal{G}_p=Gal(\Qpbar/\Qp)$. $\mathcal{G}_p$ acts on $T_p(G^0)$ via the character $\lambda(U_p^{-1})\cdot\chi_p$ (where $\chi_p$ is the $p$-adic character of $\mu_{p^\infty}$) and on $T_p(G^e)$ via the character $\lambda(U_p\cdot\dia{p}_N^{-1})\cdot\chi^{2-k}$ (where this $\chi$ really means $\chi_{{}_T}\circ\chi$).\\
3. From 1) we have an exact sequence of $\Fp$-vector space schemes $0\rightarrow G^0[p]\rightarrow G[p]\rightarrow G^e[p]\rightarrow0$ over $\Qp$ with flat extensions to $\Zp[\zeta_p]$.\\
4. $G^0[p]$ and $G^e[p]$ both have dimension $1$ and the action of $\mathcal{G}_p$ is given by $\lambda(1/a_p)\cdot\chi$ and $\lambda(a_p/\epsilon(p))\cdot\chi^{2-k}$ respectively.
\end{theorem}
\begin{proof}
Gross proves that $G$ has good reduction over $\Zp[\zeta_p]$.
It then follows from the theory of $p$-divisible groups \cite{Ta} that we have the exact sequences in 1 and 3.
He also proves that the action of Galois is as given, that the height of $G$ is $2$, and that the dimensions of $G^0[p]$ and $G^e[p]$ are each at least $1$.
But by definition the $p$-torsion of any $p$-divisible group of height $h$ always has order $p^h$.
So we see immediately that $G[p]$ has dimension $2$ and consequently that the connected and \'{e}tale components each have dimension exactly $1$.
By applying the same reasoning in reverse we see that $G^0$ and $G^e$ must then have height $1$.
\end{proof}
\begin{note} The reason why one can go back and forth between heights and dimensions as we have is that $G$ {\em always} has a $p$-divisible structure.
To be more precise, $T_p(G)$ is always free over $\Zp$.
More generally, one would like to know that $T_p(G)$ is free over $H_m$.
Unfortunately $G$ does not always have an $m$-divisible structure, although it does in many more cases than simply when $H_m=\Zp$.
\end{note} 

Now, we want to relate the splitting of $\rho_{f,p}$ to the extension class of the exact sequence
$$0\rightarrow G^0[p]\rightarrow G[p]\rightarrow G^e[p]\rightarrow0$$
so first we need to relate the representation on $G[p]$ with $\rho_f$.
$G[p]$ is precisely the $\Fp$-subspace of $J_1(Np)[p]$ on which Hecke acts according to the eigenvalues and character of $f$.
Indeed, this was the action which defined the maximal ideal $m$.
This was how we also defined the subspace $V_f$ so one might expect that the two representations are the same.
However, we had enlarged Hecke slightly.
Therefore, $G[p]$ is in fact a $2$-dimensional subspace of $V_f$ fixed by Galois.
Now, assume $\rho_f$ is irreducible so that the semisimplification of $V_f\otimes\epsilon\chi^{k-2}$ is the direct sum of copies of $\rho_f$.
Then we must have that $G[p]$ is also irreducible and furthermore that the representation on $G[p]$ is simply $\rho_f\otimes(\epsilon\chi^{k-2})^{-1}$.
This leads us to the following theorem relating the splitting of $\rho_{f,p}$ and the splitting of the exact sequence involving $G[p]$.

\begin{theorem}\label{rhosplit}Suppose $H_m=\Zp$ and that we are in the exceptional case.
Then the following are equivalent:\\
1) The canonical exact sequence involving $G[p]$ is split.\\
2) $G[p]$ is split.\\
3) $\mathcal{G}_p$ acts on $G[p]$ via the scalar $\lambda(1/a_p)\cdot\chi$.\\
4) $\rho_{f,p}$ is split.\\
5) $\rho_{f,p}$ is the scalar $\lambda(a_p)$.\\
\end{theorem}
\begin{proof}
We have already seen that $G[p]$ and $\rho_{f}$ are actually twists of each other by the character $(\epsilon\chi^{k-2})^{-1}$, which restricts to $\lambda(1/a_p^2)\cdot\chi$ in the exceptional case.
So going back and forth between the two representations is trivial.
The nontrivial part of this theorem is showing that $G[p]$ is split, meaning that there is {\em some} split exact sequence involving $G[p]$, if and only if the {\em particular} exact sequence under consideration is split.

Essentially this follows from the fact that the characters on $G^0[p]$ and $G^e[p]$ are not distinct.
To be more precise, from Theorem \ref{Gstr} we know that $\mathcal{G}_p$ acts on $G^0[p]$ and $G^e[p]$ via the characters $\lambda(1/a_p)\cdot\chi$ and $\lambda(a_p/\epsilon(p))\cdot\chi^{2-k}$.
These characters are actually identical in the exceptional case since $\epsilon(p)=a_p^2$ and $\chi^{2-k}=\chi^{1-(p-1)}=\chi$.
So $\mathcal{G}_p$ acts on the semisimplification of $G[p]$ via the scalar $\lambda(1/a_p)\chi$.
This means that if any exact sequence extending $G[p]$ is split, the representation is in fact that scalar, which in turn implies that every exact sequence is split.
\end{proof}
\section{Splitting Invariants $q_p$ and $q$}
$G^0[p]$ is a connected one dimensional $\Fp$-vector space scheme and $G^e[p]$ an \'{e}tale one dimensional $\Fp$-vector space scheme.
Therefore by Theorem \ref{Gstr} they are simply twists of $\mu_p$ and $\Z/p\Z$ by the characters $\lambda(1/a_p)$ and $\lambda(1/a_p)\cdot\chi$ respectively.
If we base extend up to a suitable field, this twisting becomes trivial and we can characterize the splitting of the exact sequence by looking at the invariant of the sequence in the group $Ext^1(\Z/p\Z,\mu_p)$.
We must choose the extension field to be large enough to trivialize the twisting and yet hopefully small enough so that the splitting status does not change.
Fortunately there is such a field $L_0$ which we will now describe.

Let $n$ be the order of $a_p$ in $\Fp^*$ and consider the following homomorphisms of groups:
$$Gal(\Qpbar/\Qp(\zeta_p))\rightarrow Gal(\Fpbar/\Fp)\rightarrow \Z/n\Z$$
The first map is the usual reduction map and the second is simply modding out by $\phi^n$ where $\phi$ is the Frobenius automorphism.
Let $\Delta_n$ be the kernel of the composition and define $L_0$ to be the fixed field of $\Delta_n$.
Since $\Delta_n$ is a normal subgroup, $L_0$ is a normal extension of $\Qp(\zeta_p)$ and we have $Gal(L_0/\Qp(\zeta_p))\cong\Z/n\Z$.
Also, since $\Delta_n$ contains the kernel of reduction, $L_0/\Qp(\zeta_p)$ is unramified.
Finally, while $\chi$ becomes trivial upon extending up to $\Qp(\zeta_p)$, the unramified character $\lambda(1/a_p)$ becomes trivial upon extending up to $L_0$ since it maps $\phi^n\rightarrow a_p^n=1$.
Now, we want to show that the splitting of the exact sequence remains invariant under this base extension, ie. the following theorem.

\begin{theorem} The exact sequence of vector space schemes is split over $\Qp$ iff it is split over $L_0$.
\end{theorem}
\begin{proof}
To prove this we first note that $L_0/\Qp$ is a finite Galois extension of degree $n(p-1)$.
Therefore for any element $\gamma\in Gal(\Qpbar/\Qp)$, $\gamma^{n(p-1)}$ must fix $L_0$.
Now suppose the sequence is split over $L_0$ and that $\gamma$ acts on $G[p]$ via the matrix
$$\left[\begin{matrix}
         a&b\\
         0&a
\end{matrix}\right]$$
with respect to any basis compatible with the exact sequence of vector space schemes.
Recall that the characters of $G^0[p]$ and $G^e[p]$ were equal in the exceptional case, so $\gamma$ must act via such a matrix.
Then $\gamma^{n(p-1)}$ must act via the matrix:
$$\left[\begin{matrix}
         a^{n(p-1)}&n(p-1)a^{n(p-1)-1}b\\
         0&a^{n(p-1)}
\end{matrix}\right]$$
But this matrix must be diagonal since $\gamma^{n(p-1)}$ fixes $L_0$ over which the sequence is split (and the representation is a scalar).
Since $n(p-1)$ is prime to $p$ this implies $b=0$ and the theorem is proved.
\end{proof}
Since $\chi$ and $\lambda(1/a_p)$ are both trivial over $L_0$, we have isomorphisms $\alpha:G^0[p]\rightarrow\mu_p$ and $\beta:\Z/p\Z\rightarrow G^e[p]$ over $L_0$.
Therefore via $\alpha$ and $\beta$ the exact sequence of vector space schemes defines a class in $Ext^1_{L_0}(\Z/p\Z,\mu_p)=L_0^*/{L_0^*}^p$.
But the vector space schemes actually had flat extensions over $\Zp[\zeta_p]$.
Therefore if we let $R_0\subset L_0$ be the ring of integers, we can extend $\alpha$ and $\beta$ to $R_0$.
This gives us a class $q_p(\alpha,\beta)\in Ext^1_{R_0}(\Z/p\Z,\mu_p)=R_0^*/{R_0^*}^p$.
Since any element of $R_0$ which is a $p$th power in $L_0$ must be a $p$th power in $R_0$, we see that we could add a third equivalent condition to the theorem.
\begin{theorem}\label{qpsplit}The following are equivalent:\\
1) The exact sequence is split over $\Qp$\\
2) The exact sequence is split over $L_0$\\
3) The exact sequence is split over $R_0$, ie. $q_p(\alpha,\beta)$ is trivial.
\end{theorem}

It is interesting and useful to note what values in particular are possible for $q_p$.
To do this we first consider how $Gal(L_0/\Qp)$ acts on the various vector space schemes.
On $G^0[p]$ and $G^e[p]$ the action is given by $\lambda(1/a_p)\cdot\chi$, and on $\mu_p$ and $\Z/p\Z$ the action is given by $\chi$ and the identity.
Therefore on $\alpha$ and $\beta$ the actions are $\lambda(a_p)$ and $\lambda(1/a_p)\cdot\chi$.
Since push-out and pull-back commute with scalar multiplication, Galois acts on the class $q_p\in Ext^1_{R_0}(\Z/p\Z,\mu_p)$ via the product $\lambda(a_p)\lambda(1/a_p)\cdot\chi=\chi$.
In other words, the values of $q_p$ must lie in the $\chi$-eigenspace.
This gives us a very useful starting point when we attempt to calculate $q_p$ and prove the main theorem.

So far the focus has been on the exact sequence of vector space schemes.
It is possible to do the analogous construction with the exact sequence of $p$-divisible groups.
In particular, $G^0$ and $G^e$ become simply $\mu_{p^\infty}$ and $\Qp/\Zp$ upon base extension to the completion of the maximal unramified extension of $\Qp(\zeta_p)$.
Let $L$ denote this field, and $R$ its ring of integers.
Then for any isomorphisms $\alpha:G^0\rightarrow\mu_{p^\infty}$ and $\beta:\Qp/\Zp\rightarrow G^e$ over $R$ we get a class $q\in Ext^1_R(\Qp/\Zp,\mu_{p^\infty})=1+\pi R$, where we choose $\pi$ to be the uniformizer $1-\zeta_p$.
Furthermore, if these isomorphisms are chosen to be compatible with the isomorphisms of vector space schemes in the obvious sense, the reduction of $q$ (mod ${R^*}^p$) is simply the image of $q_p$ under the map $R_0^*/{R_0^*}^p\rightarrow R^*/{R^*}^p$.
Therefore by the same argument this class must still be in the $\chi$-eigenspace.
However, calculating $q_p$ over $R$ would not tell us about the splitting of the exact sequence of vector space schemes over $R_0$.
It is possible, though, to calculate $q_p$ over $R_0$ from $q$, using a result of Coleman.
This will be a key element in the proof of the main theorem, along with the two inner product formulas for approximating $q$.

\section{Computing $q_p$ from $q$}

It is natural to ask what relationship exists between the splitting invariants $q$ and $q_p$ attached to the $p$-divisible group $G$.
As was already pointed out, it is clear that the images of $q$ and $q_p$ are equal in $R^*/{R^*}^p$ under the obvious maps whenever the isomorhisms $\alpha$ and $\beta$ are chosen compatibly.
But $R_0^*/{R_0^*}^p\rightarrow R^*/{R^*}^p$ is far from an injection.
So it would be impossible to determine $q_p$ from $q$ using this fact alone.
However, the two invariants are much more closely related, and in fact we {\em can} calculate $q_p$ from $q$.
The following theorem does exactly that.
Although originally proven by Coleman, here we follow the general line of reasoning of a proof by De Shalit.
\begin{theorem}\label{qpfromq} Let $\gamma:Gal(L/L_0)\rightarrow \Zp^*$ be the character $\lambda(U_p^2/\dia{p}_N)$. Choose $w\in R^*$ satisfying $w^{\sigma-1}=q^{(\gamma(\sigma)-1)/p}\quad\forall\sigma\in Gal(L/L_0)$. Then $q/w^p$ is in $R_0^*$ and $q_p\equiv q/w^p$ in $R_0^*/{R_0^*}^p$.
\end{theorem}

Before proving the theorem we should first note why there is such a $w$.
From \cite[Prop. 14.4]{Gr} we know that $q^\sigma=q^{\gamma(\sigma)}$, and we know that $p|\gamma(\sigma)-1$ because $\epsilon(p)=a_p^2$ in the exceptional case.
Therefore we have the equation
$$q^{(\gamma(\sigma)-1)/p}\left({q^{(\gamma(\tau)-1)/p)}}\right)^\sigma=q^{(\gamma(\sigma)-1)/p}q^{(\gamma(\sigma)\gamma(\tau)-\gamma(\sigma))/p}=q^{(\gamma(\sigma\tau)-1)/p}$$
This shows that $q^{(\gamma(\sigma)-1)/p}$ is a cocycle of $Gal(L/L_0)$ acting on $L^*$.
Therefore it must be a coboundary, which implies that $w$ exists and can in fact be taken to be a unit.
Furthermore, raising $w^\sigma/w$ to the $p$th power we see that
$$\left(\frac{w^\sigma}{w}\right)^p=\left(q^{(\gamma(\sigma)-1)/p}\right)^p\Rightarrow\frac{(w^p)^\sigma}{w^p}=\frac{q^\sigma}{q}\Rightarrow\left(\frac{q}{w^p}\right)^\sigma=\frac{q}{w^p}$$
But this means $q/w^p$ is in $L_0$ and hence in $R_0^*$.
Now we will show that $q/w^p$ does in fact reduce to $q_p$ in $R_0^*/{R_0^*}^p$.
\begin{proof}Recall that we have an exact sequence of $p$-divisible groups over $R_0$, namely
$$0\rightarrow G^0\rightarrow G\rightarrow G^e\rightarrow0$$
which gives us via $\alpha$ and $\beta$ a class in \mbox{$Ext^1_{R_0}(\Qp/\Zp(\lambda(U_p/\dia{p}_N)),\mu_{p^\infty}(\lambda(U_p^{-1})))$}.
By base extending to $R$ we obtain the class $q$ and by restricting to the $p$-torsion we obtain the class $q_p$.
Twisting everything by the character $\lambda(\dia{p}_N/U_p)$ of $Gal(L/L_0)$ changes neither the $p$-torsion (by choice of $L_0$) nor the base extension.
A small advantage of doing this is that $\lambda(U_p^{-1})\lambda(\dia{p}_N/U_p)=\gamma^{-1}$.
However, there is also the greater advantage that $Ext^1_{R_0}(\Qp/\Zp,\mu_{p^\infty}\gamma^{-1})=H^1(R_0,\mu_{p^\infty}\gamma^{-1})$, so we may now phrase the problem in the language of cohomology.
We have a class $\eta=\eta_q\in H^1(R_0,\mu_{p^\infty}\gamma^{-1})$ which becomes (by $res$) $q$ in $H^1(R,\mu_{p^\infty})=1+\pi R$ and reduces to $q_p$ in $H^1(R_0,\mu_p)=R_0^*/{R_0^*}^p$.
Because $\mu_{p^\infty}(L)$ is finite, $res$ is an injection.
Therefore, all we need to do is find a class which restricts to $q$ (which must be $\eta$) and show that it reduces to $q/w^p$.

To define $\eta$ we first pick a compatible system of roots of $q$ in $\Qpbar$.
By compatible, we simply mean that
$$\left(q^{1/p^{k+1}}\right)^p=q^{1/{p^k}}$$
In $H^1(R_0,\mu_{p^k}\gamma^{-1})$ we define
$$\eta_k(\sigma)=\frac{\left(q^{1/{p^k}}\right)^{\sigma\gamma^{-1}(\sigma)}}{q^{1/{p^k}}}$$
Since $q$ is fixed by $\sigma\gamma^{-1}(\sigma)$, we know $\eta_k$ must be a $p^k$th root of unity.
To show $\eta_k$ is a cocycle (of $\mu_{p^k}\gamma^{-1}$) we simply calculate
$$\eta_k(\sigma)\eta_k(\tau)^\sigma=\frac{\left(q^{1/p^k}\right)^{\sigma\gamma^{-1}(\sigma)}}{q^{1/p^k}}\cdot\frac{\left(q^{1/p^k}\right)^{\sigma\gamma^{-1}(\sigma)\tau\gamma^{-1}(\tau)}}{\left(q^{1/p^k}\right)^{\sigma\gamma^{-1}(\sigma)}}=\frac{\left(q^{1/p^k}\right)^{\sigma\tau\gamma^{-1}(\sigma\tau)}}{q^{1/p^k}}=\eta_k(\sigma\tau)$$
By the compatibility condition, the map $H^1(R_0,\mu_{p^{k+1}}\gamma^{-1})\xrightarrow{p}H^1(R_0,\mu_{p^k}\gamma^{-1})$ takes $\eta_{k+1}$ to $\eta_k$.
Therefore by taking the inverse limit we obtain a class $\eta\in H^1(R_0,\mu_{p^\infty}\gamma^{-1})$.
It is important to note that a different compatible family of roots of $q$ would define a different cocycle $\eta'$ but that $\eta/\eta'$ is actually a coboundary determined by a compatible system of roots of unity.
Therefore the class of $\eta$ is uniquely determined.

The only task remaining is to calculate the image of $\eta$ in $H^1(R,\mu_{p^\infty})$ under $res$ and the class of $\eta_1$ in $H^1(R_0,\mu_p)$.
Since $\gamma$ is trivial over $L$, we have
$$res(\eta)(\sigma)=\varprojlim\frac{\left(q^{1/p^k}\right)^\sigma}{q^{1/p^k}}$$
But this is precisely the class corresponding to $q$ in $1+\pi R$.
On the other hand, since $\eta_1(\sigma)$ is a $p$th root of unity and hence is fixed by $\gamma$, in $H^1(R_0,\mu_p)$ we have
$$\eta_1(\sigma)=\frac{(q^{1/p})^{\sigma\gamma^{-1}(\sigma)}}{q^{1/p}}=\frac{(q^{1/p})^\sigma}{(q^{1/p})^{\gamma(\sigma)}}=\frac{(q^{1/p})^\sigma}{q^{1/p}}\cdot q^{(1-\gamma(\sigma))/p}=$$
$$\frac{(q^{1/p})^\sigma}{q^{1/p}}\cdot\frac{w}{w^\sigma}=\frac{\left(\frac{q^{1/p}}{w}\right)^\sigma}{\left(\frac{q^{1/p}}{w}\right)}=\frac{((q/w^p)^{1/p})^\sigma}{(q/w^p)^{1/p}}$$
But this is precisely the class corresponding to $q/w^p$ in $H^1(R_0,\mu_p)=R_0^*/{R_0^*}^p$.
Therefore we have proved the theorem.
\end{proof}
\section{Formulas for $\log q$ and $d\log q$}
At this point we know that we could (in theory) determine whether $\rho_{f,p}$ is split by choosing isomorphisms $\alpha$ and $\beta$ and then determining the splitting invariant $q(\alpha,\beta)\in 1+\pi R$.
We still do not have a practical criterion, however, which can be checked by straightforward calculations.
In this section we will introduce two formulas which {\em can} be evaluated by straightforward calculations.
The first formula, taken directly from \cite[Thm. 4.4]{CV}, uses the cup product on $H^1_{DR}(I)$ (where $I=I_1(N)$ is the Igusa curve of level $N$ in characteristic $p$) to compute $d\log q$. 
The second formula, taken directly from \cite[Thm. 2.1]{C1}, uses an inner product defined on a subspace of $H^0(X_1(Np),\Omega)$ to compute $\log q$.
The main theorem will then be a statement of how the splitting of $\rho_{f,p}$ is precisely related to the triviality of these inner products.

The first step to understanding the formulas is to give a different but equivalent interpretation of $\alpha$ and $\beta$.
The choice of any homomorphism from $\Qp/\Zp$ to $G^e$ over $R$ is equivalent to choosing an element of $T_pG^e=T_p{\bar G}$.
When $H_m=\Zp$, we have seen that this is a free $\Zp$-module of rank $1$.
Therefore, $\beta$ is an isomorphism exactly when it corresponds to a generator of $T_p{\bar G}$.
Similarly, a choice of homomorphism from $G^0$ to $\mu_{p^\infty}$ over $R$ is equivalent to choosing an element of $T_p{\bar G}'$, where $G'=Hom(G,\mu_{p^\infty})$ is the Cartier dual of $G$.
Again, since this is a free $\Zp$-module of rank $1$, an isomorphism $\alpha$ corresponds to a generator of $T_p{\bar G}'$.
With this interpretation of $\alpha$ and $\beta$, $q(\alpha,\beta)$ is just the usual Serre-Tate invariant.
$$q:T_p{\bar G}\times T_p{\bar G}'\rightarrow 1+\pi R$$
If we follow this by the $p$-adic logarithm and extend by scalars to $R$, this is what is meant by the map
$$\log q:(T_p{\bar G}\otimes_{\Zp} R)\times(T_p{\bar G}'\otimes_{\Zp}R)\rightarrow R$$
If instead we follow it by the map $d\log:a\rightarrow da/a$ and extend by scalars, this is what is meant by the map
$$d\log q:(T_p{\bar G}\otimes_{\Zp} R)\times(T_p{\bar G}'\otimes_{\Zp}R)\rightarrow \Omega_{R/\Zp^{\text{unr}}}$$
Here $\Zp^{\text{unr}}$ is the ring of integers in the completion of the maximal unramified extension of $\Qp$.

Now, one has to be a little careful in that the $p$-divisible group $G=G_f$ which we have defined and the $p$-divisible group referred to as $G$ in \cite[Thm. 4.4]{CV} and \cite[Thm. 2.1]{C1} are not the same.
While our $G$ was defined by $T_pG=\epsilon_m T_p(J_1(Np))$, the other is defined by
          $$T_pG=\bigcap_{n}U_p^n \left(T_p(J_1(Np))^Z\right)$$
By the definition of $m$, though, the operator $U_p\equiv a_p$ (mod $m$) on $G_f$, which means in particular $U_p$ is invertible.
Also, by the proof of \cite[Prop. 12.9]{Gr}, the points of $G_f$ are indeed killed by the correspondence $Z=\sum_{d\in\Fp^*}\dia{d}_p$.
So {\em our} $G=G_f$ is simply the subgroup of the larger $G$ which is cut out by the idempotent $\epsilon_m$.
In \cite[Thm. 4.4]{CV} and \cite[Thm. 2.1]{C1} it is also stated that 
         $$T_pG'=\bigcap_{n}{U_p'}^n\left(T_p(J_1(Np))^Z\right)$$
Although this is not stated explicitly, it is clear that the canonical pairing of $T_pG$ and $T_pG'$ into $\Zp(1)$ in that case simply comes from the Weil pairing on the Tate module of $J_1(Np)$.
This means that when we apply the theorems we are implicitly also identifying our $T_pG'$ with $\epsilon_{m'}T_p(J_1(Np))$ where $m'=ros(m)$ is the image of $m$ under the Rosati involution of $End(J_1(Np))$.

The last point which is essential for understanding the two formulas is that there is a correspondence between elements of $T_p{\bar G}\otimes R$ (or $T_p{\bar G}'\otimes R$) and differentials on $X$, the canonical model for $X_1(Np)$, over $R$.
In particular, applying \cite[Lemma 4.3]{CV} to our $G$ and $G'$, we see that there are natural isomorphisms
$$T_p{\bar G}'\otimes_{\Zp}R\rightarrow \epsilon_m H^0(X,\Omega_{X/R})=\Omega_G$$
$$T_p{\bar G}\otimes_{\Zp}R\rightarrow \epsilon_{m'}H^0(X,\Omega_{X/R})=\Omega_{G'}$$
The subtlety and the power of these isomorphisms is really in the integral structures over $R$.
When $H_m=\Zp$ and the lifting $F$ of $f$ is unique, $\Omega_G$ is generated over $R$ by the regular differential $\omega_F$.
Likewise $\Omega_{G'}$ is generated by $\omega_{F|w}$ for the automorphism $w=w_{\zeta_p}$ of $X_1(Np)$ (see \cite[Prop. 8.4, 6.14]{Gr}).
Therefore, by these identifications, the original isomorphisms $\alpha$ and $\beta$ correspond up to units in $R^*$ to $\omega_F$ and $\omega_{F|w}$.
The integral structure also makes it possible to reduce a differential on $X$ (mod $\pi$) to obtain a differential on the Igusa curve $I=I_1(N)$ by \cite[Prop 7.1]{Gr}.
Aside from respecting the integral structures, these isomorphisms are also nice in that they commute with Hecke in the only possible sense, namely $h\beta=ros(h)\omega_\beta$ and $ros(h)\alpha=h\omega_\alpha$.
With these identifications in mind then, it now makes sense to state the first formula.
\begin{theorem}
If $\beta\in T_p{\bar G}\otimes_{\Zp}R$, $\alpha\in T_p{\bar G}'\otimes_{\Zp}R$, and $\omega_\alpha|_I=\omega_f$, then
$$d\log q(\alpha,\beta)=a_p<w^*\omega_\beta|_I,[f']>_I d\pi+\cdots$$
\end{theorem}
\begin{corollary}\label{dlogq}
If $\beta\in T_p{\bar G}\otimes_{\Zp}R$, $\alpha\in T_p{\bar G}'\otimes_{\Zp}R$, and $\omega_\alpha|_I=\omega_f$, then
$$d\log q\equiv 0\pmod{\pi d\pi}\ \text{iff}\ <w^*\omega_\beta|_I,[f']>_I=0$$
\end{corollary}
This is almost an exact restatement of \cite[Thm. 4.4]{CV}.
The simplifications come from the fact that we are in the exceptional case, and that our group $G$ is less general.
We should note that in the exceptional case, $f'$ refers to the form $\theta f$ on $X_1(N)$.
While $f$ actually defines a holomorphic differential on $I$ (the reduction of $\omega_F$), both $f$ and $f'$ are shown in \cite[Sect. 13]{Gr} to define classes in $H^1_{DR}(I)$.
Also, it should be noted that while our choice of uniformizer $\pi$ is not the same as the one in the theorem, the corollary is independent of that choice and is all that we will actually need.

The second formula makes use of an inner product, denoted $(\ ,\ )_\infty$ and introduced in \cite{C1}, on a subspace $M_0$ of $H^0(X_1(Np),\Omega)$.
The subspace actually contains the $W^{\text{ord}}$ and $W^{\text{anti-ord}}$ of \cite{CV} and \cite{C1} so in particular it certainly contains our $\Omega_G$ and $\Omega_{G'}$.
This time the statement which we need is a precise quote of \cite[Theorem 2.1]{C1}.
\begin{theorem}\label{logq} Let $\alpha\in T_p{\bar G}'\otimes_{\Zp}R$ and $\beta\in T_p{\bar G}\otimes_{\Zp}R$. Then
            $$(\omega_\alpha,\omega_\beta)_\infty=\log q(\alpha,\beta)$$
\end{theorem}

\section{The Main Theorem}

\begin{theorem}
Suppose $H_m=\Zp$ and we are are in the exceptional case. Then $\rho_{f,p}$ is split iff $(\omega_F,\omega_{F|w})_\infty\equiv 0\pmod{\pi^{p+1}}$ and $<[f],[f']>_I=0$.
\end{theorem}
\begin{proof}
We begin by choosing any isomorphisms 
$$\alpha:G^0\rightarrow\mu_{p^\infty}\ \text{and}\ \beta:\Qp/\Zp\rightarrow G^e$$
over $R$, compatible with a fixed choice of isomorphisms
$$\alpha:G^0[p]\rightarrow\mu_p\ \text{and}\ \beta:\Z/p\Z\rightarrow G^e[p]$$
over $R_0$. 
By Theorems \ref{rhosplit} and \ref{qpsplit} we know $\rho_{f,p}$ is split iff $q_p(\alpha,\beta)=0$.
On the other hand, Theorem \ref{qpfromq} gives us a way to calculate $q_p$ from $q$, and we have inner product formulas for $d\log q$ and $\log q$ from the previous section.
This line of reasoning gives the proof of the main theorem its overall structure.

To obtain a starting point for $q$ then, we first use the fact that the reduction of $q$ must be in the $\chi$-eigenspace of the Galois module
$$R^*/{R^*}^p=(1+\pi R)/(1+\pi R)^p=(1+\pi R)/(1+\pi^p R)$$
This last equality follows from the fact that $R$ is the ring of integers in the completion of the maximal unramified extension of $Q_p(\zeta_p)$.
Using the triviality of Frobenius and solving for the $\pi^i$ coefficients iteratively one can show that the $\chi$-eigenspace is precisely the $p$th roots of unity.
So we {\em could} take as our starting point $q=\zeta^s(1+r\pi^p)$, for some $s\in \Z/p\Z$, and $r\in R$.
It is possible to be more precise, though, about the $r$ in the expression.
Recall that $q$ is actually in the $\lambda(A_p^2/\epsilon(p))\cdot\chi$-eigenspace of the Galois module $1+\pi R$, and $A_p^2/\epsilon(p)\equiv1$.
So acting on $q$ by any Frobenius automorphism $\phi$ we have
\begin{align*}
                    q^\phi&=q\cdot\left(q^p\right)^{(A_p^2/\epsilon(p)-1)/p}\\
  (\zeta^s(1+r\pi^p))^\phi&=\zeta^s(1+r\pi^p)\left((\zeta^s(1+r\pi^p))^p\right)^{(A_p^2/\epsilon(p)-1)/p}\\
    \zeta^s(1+r^\phi\pi^p)&=\zeta^s(1+r\pi^p)((1+r\pi^p)^p)^{(A_p^2/\epsilon(p)-1)/p}\\
       \zeta^s(1+r^p\pi^p)&\equiv\zeta^s(1+r\pi^p)\quad(\text{mod}\ \pi^{p+1})
\end{align*}
But this means that the reduction of $r$ (mod $\pi$) is actually in $\Fp$.
So we may actually start with the congruence
$$q\equiv\zeta^s(1+\pi^p)^t\pmod{\pi^{p+1}}$$
with both $s$ and $t$ in $\Z/p\Z$.

To calculate $q_p$ from $q$ using Theorem \ref{qpfromq} we need to find a $w\in R^*$ satisfying $w^{\sigma-1}=q^{(\gamma(\sigma)-1)/p}$ for all $\sigma\in Gal(L/L_0)$.
Since this Galois group is generated by a power of $\phi$, it suffices to find $w$ such that $w^{\phi-1}=q^{(A_p^2/\epsilon(p)-1)/p}$.
The map $\phi-1$ from $1+\pi R$ to itself is surjective, so we can do this in pieces by choosing $u,v\in 1+\pi R$ which satisfy 
$$u^{\phi-1}=\zeta=1-\pi\qquad\qquad v^{\phi-1}=1+\pi^p$$
Then we can simply let $w=(u^s v^t)^{(A_p^2/\epsilon(p)-1)/p}$.
Immediately we see that 
\begin{align*}
   u=1+u_0\pi\ &\Rightarrow\ u_0^p\equiv u_0-1\pmod{\pi}\\
   v=1+v_0\pi\ &\Rightarrow\ v_0^p\equiv v_0\pmod{\pi}
\end{align*}
Since $(1+r\pi)^p\equiv1+(r^p-r)\pi^p\pmod{\pi^{p+1}}$ for any $r$, this implies
\begin{align*}
u^p&\equiv(1-\pi^p)\equiv(1+\pi^p)^{-1}\pmod{\pi^{p+1}}\\
v^p&\equiv1\pmod{\pi^{p+1}}
\end{align*}
Plugging in for $w^p$ and then for $q$ using Theorem \ref{qpfromq} we get
\begin{align*}
w^p&\equiv(1+\pi^p)^{-s(A_p^2/\epsilon(p)-1)/p}\pmod{\pi^{p+1}}\\
q_p&\equiv\zeta^s(1+\pi^p)^{t+s(A_p^2/\epsilon(p)-1)/p}\pmod{\pi^{p+1}}
\end{align*}

\begin{lemma}$\rho_{f,p}$ is split iff $s=t=0$.
\end{lemma}
\begin{proof}[Proof (of lemma)]
If $s=t=0$, $q_p$ is in $1+\pi^{p+1}R$ and is therefore a $p$th power.
By Theorems \ref{rhosplit} and \ref{qpsplit} this means $\rho_{f,p}$ is split.
Conversely, we know that $(1+\pi R)^p\subset(1+\pi^p R)$, so if $q_p$ is a $p$th power it follows immediately that $s=0$.
Furthermore, it is easily shown that a $p$th root of $1+t\pi^p+\cdots$ for $t\in\Z/p\Z$ generates a degree $p$ extension of $\Qp(\zeta_p)$ unless $t=0$.
But the degree of $L_0/\Qp(\zeta_p)$ was prime to $p$.
Therefore, once $s=0$, it is clear that $t$ must also be $0$.
\renewcommand{\qedsymbol}{}
\end{proof}
\begin{lemma}$s=t=0$ iff $d\log q\equiv0\pmod{\pi d\pi}$ and $\log q\equiv0\pmod{\pi^{p+1}}$.
\end{lemma}
\begin{proof}[Proof (of lemma)]
Regardless of $t$, it is always true that 
$$d\log q\equiv(-s/\zeta)d\pi\pmod{\pi d\pi}$$
While the $\log q$ congruence is not independent of $s$ in general, in the case that $s=0$ we have simply 
$$\log q\equiv t\pi^p\pmod{\pi^{p+1}}$$
This is actually all we need to know logically to prove the lemma.
First, suppose $s=0$ and $t=0$.
The above calculations clearly imply the lemma in that case.
Conversely, suppose that $d\log q$ and $\log q$ are sufficiently trivial.
It follows immediately from the first calculation that $s=0$.
This in turn validates the second calculation which implies $t=0$ and we are done.
\renewcommand{\qedsymbol}{}
\end{proof}
The proof of the main theorem now comes down to a relatively simple application of the two inner product formulas.
Since $\alpha$ and $\beta$ correspond to generators of $T_p{\bar G}'$ and $T_p{\bar G}$ respectively, we must have
$$\omega_\alpha=r_1\omega_F,\quad\omega_\beta=r_2\omega_{F|w}$$
for units $r_1,r_2\in R^*$.
Since multiplication by units does not affect the triviality of either inner product we have by Corollary \ref{dlogq} and Theorem \ref{logq}
$$d\log q\equiv0\pmod{\pi d\pi}\ \Leftrightarrow\ <[f],[f']>_I$$
$$\log q\equiv0\pmod{\pi^{p+1}}\ \Leftrightarrow\ (\omega_F,\omega_{F|w})_\infty\equiv0\pmod{\pi^{p+1}}$$
Combining these last equivalences with the two lemmas, the Main Theorem follows.
\end{proof}
\begin{remark}
In \cite{Mc} these inner products are explicitly calculated for two examples.
The first begins with a weight $5$ form for $\Gamma_1(4)$ and the second with a weight $7$ form for $\Gamma_1(3)$.
Unfortunately, however, both associated representations are in fact reducible and come from companion forms which are weight $1$ Eisenstein series.
So while the calculations illustrate methods for computing the splitting criterion, the necessary hypotheses of the criterion are not satisfied.
\end{remark}
\noindent
{\em Acknowledgements.} The author would like to thank Robert Coleman and Benedict Gross whose cited works provided both essential results and helpful exposition of the subject.
To Robert Coleman, especially, we are deeply indebted for countless helpful conversations and truly exceptional advising.
.


\begin{thebibliography}{99}
\bibitem[C1]{C1}
R. Coleman,
\emph{A $p$-Adic Inner Product on Elliptic Modular Forms},
Proceedings of the Barsotti Conference,
pp. ~125--151.
\bibitem[C2]{C2}
R. Coleman,
\emph{Reciprocity Laws on Curves},
Compositio Mathematica \textbf{72}~(1989), 205--235.
\bibitem[CV]{CV}
R. Coleman, J. Voloch,
\emph{Companion Forms and Kodaira-Spencer Theory},
Invent. Math. \textbf{110}~(1992), 263--281.
\bibitem[De]{De}
P. Deligne,
\emph{Formes Modulaires et representations $l$-adiques},
Seminaire Bourbaki \textbf{355},
Lecture Notes in Math. \textbf{179}~(1971), 139--172.
\bibitem[DR]{DR}
P. Deligne, M. Rapoport,
\emph{Schemas de Modules de Courbes Elliptiques},
Lecture Notes in Math. \textbf{349}~(1973), 143--316.
\bibitem[DI]{DI}
F. Diamond, J. Im,
\emph{Modular Forms and Modular Curves},
Canadian Math. Society Conference Proceedings,
pp. ~1-95.
\bibitem[Gr]{Gr}
B. Gross,
\emph{A Tameness Criterion for Galois Representations Associated to Modular Forms (mod $p$)},
Duke Math. Journal \textbf{61}~(1990), 445--516.
\bibitem[Ka]{Ka}
N. Katz,
\emph{$p$-Adic Properties of Modular Schemes and Modular Forms},
Lecture Notes in Math. \textbf{350}~(1973), 69--170.
\bibitem[Mc]{Mc}
K. McMurdy,
\emph{A Splitting Criterion for Galois Representations Associated to Exceptional Modular Forms},
Ph.D. thesis, University of California, Berkeley, 2001.
\bibitem[Sh]{Sh}
G. Shimura,
\emph{Introduction to the Arithmetic Theory of Automorphic Functions},
Iwanami Shoten and Princeton Univ. Press,
Princeton, 1971.
\bibitem[Ta]{Ta}
J. Tate,
\emph{$p$-Divisible Groups},
Proceedings of a Conference on Local Fields,
pp. 158--183.
\end{thebibliography}
\end{document}